\documentclass[11pt, a4paper]{article}
\usepackage{times}
\usepackage{a4wide}
\usepackage[dvips, hyperindex]{hyperref}
\usepackage[british]{babel}
\usepackage{enumerate, longtable}
\usepackage{amsmath, amscd, amsfonts, amssymb, latexsym, theorem, comment, stmaryrd}
\usepackage{xypic}
\usepackage[T1]{fontenc}
\usepackage[latin1]{inputenc}

{\theorembodyfont{\itshape}   \newtheorem{thm}{Theorem}[section]}
{\theorembodyfont{\itshape}   }
{\theorembodyfont{\itshape}   }
{\theorembodyfont{\itshape}   \newtheorem{rem}[thm]{Remark}}
{\theorembodyfont{\itshape}   }
{\theorembodyfont{\itshape}   \newtheorem{prop}[thm]{Proposition}}
{\theorembodyfont{\itshape}   \newtheorem{cor}[thm]{Corollary}}
{\theorembodyfont{\itshape}   }
{\theorembodyfont{\itshape}   }
{\theorembodyfont{\itshape}   \newtheorem{sit}[thm]{Situation}}
{\theorembodyfont{\itshape}   }
{\theorembodyfont{\itshape}   }

\newcommand{\CC}{\mathbb{C}}
\newcommand{\QQ}{\mathbb{Q}}
\newcommand{\FF}{\mathbb{F}}
\newcommand{\ZZ}{\mathbb{Z}}

\newcommand{\TT}{\mathbb{T}}

\newcommand{\NN}{\mathbb{N}}

\newcommand{\Qbar}{{\overline{\QQ}}}
\newcommand{\Fbar}{{\overline{\FF}}}

\newcommand{\Hom}{{\rm Hom}}

\newcommand{\GL}{\mathrm{GL}}
\newcommand{\Aut}{\mathrm{Aut}}

\newcommand{\grs}{\mathrm{gr.sch.}}

\DeclareMathOperator{\Gal}{Gal}

\DeclareMathOperator{\Frob}{Frob}
\DeclareMathOperator{\Tr}{Tr}
\DeclareMathOperator{\Det}{Det}

\DeclareMathOperator{\coker}{coker}

\DeclareMathOperator{\Spec}{Spec}

\newcommand{\pf}{{\bf Proof. }}
\newcommand{\qed}{\hspace* {.5cm} \hfill $\Box$}

\newcommand{\cO}{\mathcal{O}}

\newcommand{\calG}{\mathcal{G}}

\newcommand{\rk}{{\rm rk}}
\newcommand{\nr}{{\rm nr}}
\newcommand{\ev}{{\rm ev}}
\newcommand{\gl}{{\rm gl}}

\newcommand{\Halg}{{\rm HA}}

\newcommand{\m}{\mathfrak{m}}
\newcommand{\mbar}{{\overline{\mathfrak{m}}}}

\title{Multiplicities of Galois representations of weight one\\
(with an appendix by Niko Naumann)}
\author{Gabor Wiese\footnote{
NWF 1-Mathematik,
Universität Regensburg,
D-93040 Regensburg,
Germany, gabor@pratum.net}}

\begin{document}

\maketitle

\begin{abstract}
  In this article we consider mod~$p$ modular Galois representations
  which are unramified at~$p$ such that the Frobenius element at~$p$
  acts through a scalar matrix. The principal result states that the
  multiplicity of any such representation is bigger than~$1$.

\end{abstract}

\section{Introduction}

A continuous odd irreducible Galois representation $\rho:
\Gal(\Qbar/\QQ) \to \GL_2(\Fbar_p)$ is said to be of weight one if it
is unramified at~$p$. According to Serre's conjecture (with the
minimal weight as defined in~\cite{edixhoven-serre}), all such
representations should arise from Katz modular forms of weight~$1$
over~$\Fbar_p$ for the group $\Gamma_1(N)$ with $N$ the (prime to~$p$)
conductor of~$\rho$. Assuming the modularity of~$\rho$, this is known
if $p>2$ or if $p=2$ and the restriction of~$\rho$ to a decomposition
group at~$2$ is not an extension of twice the same character. A weight~$1$
Katz modular form over~$\Fbar_p$ can be embedded into weight~$p$
and the same level in two different ways: by multiplication by the
Hasse invariant of weight~$p-1$ and by applying the Frobenius
(see~\cite{edixhoven}, Section~4). Hence, the corresponding
eigenform(s) in weight~$p$ should be considered as old forms; they lie
in the ordinary part.

A modular Galois representation $\rho: \Gal(\Qbar/\QQ) \to
\GL_2(\Fbar_p)$ of conductor~$N$ can be realised with a certain
multiplicity (see Proposition~\ref{blr}) on the $p$-torsion of
$J_1(Np)$ or $J_1(N)$ (for $p=2$).  In this article we prove that this
multiplicity is bigger than~$1$ if~$\rho$ is of weight one
and~$\Frob_p$ acts by scalars. If $p=2$, we also assume that the
corresponding weight~$1$ form exists. Together with
\cite{buzzard_appendix}, Theorem~6.1, this completely settles the
question of multiplicity one for modular Galois representations. Its
study had been started by Mazur and continued among others by Ribet,
Gross, Edixhoven and Buzzard. The first example of a modular Galois
representation of not satisfying multiplicity one was found by Kilford
in~\cite{kilford-nongorenstein}. See~\cite{kilford_wiese} for a more
detailled exposition.

A systematic computational study of the multiplicity of Galois
representations of weight one has been carried out
in~\cite{kilford_wiese}. The data gathered suggest that the
multiplicity always seems to be~$2$ if it is not~$1$. Moreover, the
local factors of the Hecke algebras are becoming astonishingly large.

\subsection*{Overview}

We give a short overview over the article with an outline of the
proof. In Section~\ref{jac} an isomorphism between a certain part of
the $p$-torsion of a Jacobian of a modular curve with a local factor
of a mod~$p$ Hecke algebra is established (Proposition~\ref{propq}).
As an application one obtains a mod~$p$ version of the Eichler-Shimura
isomorphism (Corollary~\ref{gormult}). Together with a variant of a
well-known theorem by Boston, Lenstra and Ribet
(Proposition~\ref{blr}) one also gets an isomorphism between a certain
kernel in the local mod~$p$ Hecke algebra and a part of the
corresponding Galois representation. This gives for instance a precise
link between multiplicities and properties of the Hecke algebra
(Corollary~\ref{multgor}). In Section~\ref{frob} it is proved
(Theorem~\ref{thmtp}) that under the identification of
Section~\ref{jac}, the Frobenius at~$p$ on the part of the Galois
representation corresponds to the Hecke operator~$T_p$ in the Hecke
algebra. This relation is exploited in Section~\ref{conclude} to
obtain the principal result (Theorem~\ref{mainthm}) and a couple of
corollaries.

\subsection*{Notations}

For integers $N \ge 1$ and $k \ge 1$, we let $S_k(\Gamma_1(N))$ be the
$\CC$-vector space of holomorphic cusp forms and
$S_k(\Gamma_1(N),\FF_p)$ the $\FF_p$-vector space of Katz cusp forms
on~$\Gamma_1(N)$ of weight~$k$.  Whenever $S \subseteq R$ are rings,
$m$ is an integer and $M$ is an $R$-module on which the Hecke and
diamond operators act, we let $\TT_S^{(m)}(M)$ be the $S$-subalgebra
inside the $R$-endomorphism ring of~$M$ generated by the Hecke
operators $T_n$ with $(n,m)=1$ and the diamond operators. If $\phi: S
\to S'$ is a ring homomorphism, we let $\TT^{(m)}_\phi (M) :=
\TT^{(m)}_S(M) \otimes_S S'$ or with $\phi$ understood $\TT^{(m)}_{S
  \to S'} (M)$. If $m=1$, we drop the superscript.

Every maximal ideal 
$\mbar \subseteq \TT_{\ZZ \to \FF_p}(S_k(\Gamma_1(N)))$ 
corresponds to a Galois conjugacy class of
cusp forms over~$\Fbar_p$ of weight~$k$ on~$\Gamma_1(N)$.  One
can attach to~$\mbar$ by work of Shimura and Deligne a continuous odd
semi-simple Galois representation $\rho_\mbar: \Gal(\Qbar/\QQ) \to
\GL_2(\Fbar_p)$ which is unramified outside~$Np$ and satisfies
$\Tr(\rho_\mbar(\Frob_l)) \equiv T_l \mod \mbar$ and $\Det(\rho_\mbar(\Frob_l))
\equiv \langle l \rangle l^{k-1} \mod \mbar$ for all primes $l \nmid
Np$ via an embedding 
$\TT_{\ZZ \to \FF_p}(S_k(\Gamma_1(N)))/\mbar \hookrightarrow \Fbar_p$.
All Frobenius elements $\Frob_l$ are arithmetic ones.

For all the article we fix an isomorphism $\CC \cong \Qbar_p$
and a ring surjection $\overline{\ZZ}_p \to \Fbar_p$.
If $K$ is a field, we denote by $K(\epsilon) = K[\epsilon]/(\epsilon^2)$
the dual numbers. For a finite flat group scheme $G$,
the Cartier dual is denoted by~${}^t G$.
The maximal unramified extension of~$\QQ_p$ (inside $\Qbar_p$) is
denoted by~$\QQ_p^\nr$ and its integer ring by~$\ZZ_p^\nr$.

\subsection*{Situations}

We shall often assume one of the following two situations. In the
applications, the second part will be taken for~$p=2$.

\begin{sit}\label{sit}
\begin{enumerate}[(I)]

\item Let $p$ be an odd prime and $N$ a positive integer not divisible
  by~$p$. For $m \in \NN$ write $\TT^{(m)}_{\ZZ_p}$ for the Hecke
  algebra $\TT^{(m)}_{\ZZ \to \ZZ_p}(S_2(\Gamma_1(Np)))$. Let $\m$ be
  an ordinary (i.e.\ $T_p \not\in \m$) maximal ideal of $\TT_{\ZZ_p}$.
  Denote the image of~$\m$ in $\TT_{\FF_p} := \TT_{\ZZ
    \to \FF_p}(S_2(\Gamma_1(Np)))$ by~$\mbar$. Assume that $\rho_\mbar$ is
  irreducible. Let $\m^{(m)} = \m \cap \TT^{(m)}_{\ZZ_p}$ and
  similarly for $\mbar^{(m)}$.

  Let furthermore $K = \QQ_p(\zeta_p)$ and $\cO = \ZZ_p[\zeta_p]$ with
  a primitive $p$-th root of unity~$\zeta_p$.
  Also let $J := J_1(Np)_\QQ$ be the Jacobian of $X_1(Np)$ over~$\QQ$.

\item Let $p$ be any prime and $N$ a positive integer not divisible
  by~$p$. For $m \in \NN$ write $\TT_{\ZZ_p}^{(m)}$ for the Hecke
  algebra $\TT^{(m)}_{\ZZ \to \ZZ_p}(S_2(\Gamma_1(N)))$. Let $\m$ be
  an ordinary (i.e.\ $T_p \not\in \m$) maximal ideal of $\TT_{\ZZ_p}$.
  Denote the image of~$\m$ in $\TT_{\FF_p} := \TT_{\ZZ
    \to \FF_p}(S_2(\Gamma_1(N)))$ by~$\mbar$. Assume that $\rho_\mbar$ is
  irreducible. Let $\m^{(m)} = \m \cap \TT^{(m)}_{\ZZ_p}$ and
  similarly for $\mbar^{(m)}$.

  Let furthermore $K = \QQ_p$ and $\cO = \ZZ_p$.
  Also let $J := J_1(N)_\QQ$ be the Jacobian of $X_1(N)$ over~$\QQ$.

\end{enumerate}
\end{sit}

\subsection*{Acknowledgements}

I am indebted to Kevin Buzzard for his suggestions and help concerning
Section~\ref{jac} and his comments. It is a pleasure to thank Lloyd
Kilford for his work on~\cite{kilford_wiese} without which the present
article would certainly not exist. I thank Niko Naumann for useful
comments and his interest leading to his proof of Theorem~\ref{thmtp},
and Bas Edixhoven for many explanations and discussions concerning
this and related subjects.

\section{Hecke algebras, Jacobians and $p$-divisible groups}\label{jac}

Let us assume one of the two cases of Situation~\ref{sit}. The maximal
ideal~$\m$ of~$\TT_{\ZZ_p}$ corresponds to an idempotent~$e_\m \in
\TT_{\ZZ_p}$, in the sense that applying~$e_\m$ to any
$\TT_{\ZZ_p}$-module is the same as localising the module at~$\m$.
Let $\calG$ be the $p$-divisible group~$J[p^\infty]_\QQ$ over~$\QQ$.
Consider the Tate module $T_p J = T_p \calG =
\underset{\leftarrow}{\lim} J[p^n](\Qbar)$. It is a
$\TT_{\ZZ_p}[\Gal(\Qbar/\QQ)]$-module. The idempotent~$e_\m$ gives
rise to an endomorphism of $\TT_{\ZZ_p}[\Gal(\Qbar/\QQ)]$-modules
on~$T_p J$. Such an endomorphism comes from an endomorphism~$e_\m$ of
the $p$-divisible group~$\calG$, which is also an idempotent. We put
$G = e_\m \calG$ and say that this is the {\em $p$-divisible group
  over~$\QQ$ attached to~$\m$}. We shall mainly be interested in the
$p$-torsion of~$G$. However, making the detour via $p$-divisible
groups allows us to quote the following theorem by Gross.

\begin{thm}[Gross]\label{propgross}
  Assume any of the two cases of Situation~\ref{sit}. Let $G$ be the
  $p$-divisible group over~$\QQ$ attached to~$\m$, as explained above.
  Let $h = \rk_{\ZZ_p} \TT_{\ZZ_p,\m}$, where $\TT_{\ZZ_p,\m}$ denotes
  the localisation of $\TT_{\ZZ_p}$ at~$\m$.

\begin{enumerate}[(a)]

\item
  The $p$-divisible group~$G$ acquires good reduction over~$\cO$. We write
  $G_\cO$ for the corresponding $p$-divisible group over~$\cO$.
  It sits in the exact sequence
  $$ 0 \to G^0_\cO \to G_\cO \to G^e_\cO \to 0,$$
  where $G^e_\cO$ is \'etale and $G^0_\cO$ is of multiplicative type,
  i.e.\ its Cartier dual is \'etale. The exact sequence is preserved
  by the action of the Hecke correspondences.

\item
  Over $\cO[\zeta_N]$ the group $G_{\cO[\zeta_N]}$ is isomorphic to
  its Cartier dual ${}^t G_{\cO[\zeta_N]}$.
  This gives isomorphisms of $p$-divisible groups over~$\cO[\zeta_N]$
  $$ G^e_{\cO[\zeta_N]} \cong {}^t G^0_{\cO[\zeta_N]} \;\;\; \text{ and }\;\;\;
     G^0_{\cO[\zeta_N]} \cong {}^t G^e_{\cO[\zeta_N]}.$$

\item
  We have $G^e_{\Fbar_p}[p] \cong (\ZZ/p\ZZ)^h_{\Fbar_p}$ and
  $G^0_{\Fbar_p}[p] \cong \mu_{p, \Fbar_p}^h$.

\end{enumerate}
\end{thm}

\pf
(a) The statement on the good reduction is \cite{gross}, Prop.~12.8~(1)
and~12.9~(1). The exact sequence is proved in \cite{gross}, Prop.~12.8~(4)
and~12.9~(3). That it is preserved by the Hecke correspondences is a
consequence of the fact that there are no non-trivial
morphisms from a connected group scheme to an \'etale one, whence
any Hecke correspondence on~$G$ restricts to~$G^0$.

(b) The Cartier self-duality of~$G$ over~$K(\zeta_N)$ is also proved in \cite{gross}, 
Prop.~12.8~(1) and~12.9~(1). It extends to a self-duality over~$\cO[\zeta_N]$.
The second statement follows as in~(a) from the non-existence
of non-trivial morphisms from $G^0$ to~$G^e$ over~$\cO[\zeta_N]$;
this argument gives $G^0 \cong {}^t G^e$. Applying Cartier duality
to this, we also get $G^e \cong {}^t G^0$.

(c) By Part~(b), $G^e$ and~$G^0$ have equal height. That height is
equal to~$h$ by \cite{gross}, Prop.~12.8~(1) and~12.9~(1). The
statement is now due to the fact that up to isomorphism the given
group schemes are the only ones of rank~$p^h$ which are killed by~$p$
and which are \'etale or of multiplicative type, respectively.
\qed
\medskip

The last point makes the ordinarity of~$\mbar$ look like the ordinarity
of an abelian variety.

\begin{prop}\label{propq}
  Assume any of the two cases of Situation~\ref{sit} and let~$G$ be the
  $p$-divisible group attached to~$\m$. Then we have the isomorphism
  $G^0[p](\Qbar_p) \cong \TT_{\FF_p,\mbar}$ of $\TT_{\FF_p,\mbar}$-modules.
\end{prop}

\pf
Taking the $p$-torsion of the $p$-divisible groups in Thm.~\ref{propgross}~(a),
one obtains the exact sequence
\begin{equation}\label{grossdrei}
 0 \to G^0_\cO[p](\Qbar_p) \to G_\cO[p](\Qbar_p) \to G_\cO^e[p](\Qbar_p) \to 0
\end{equation}
of $\TT_{\FF_p, \mbar}$-modules with Galois action.
We also spell out the dualities in Thm.~\ref{propgross}~(b) restricted
to the $p$-torsion:
\begin{equation}\label{grossvier}
\begin{split}
G^0_{\cO[\zeta_N]}[p] & \cong \Hom_{\grs/\cO[\zeta_N]}(G^e[p]_{\cO[\zeta_N]},\mu_{p,\cO[\zeta_N]})
    \;\;\; \text{ and } \;\;\; \\
G^e_{\cO[\zeta_N]}[p] & \cong \Hom_{\grs/\cO[\zeta_N]}(G^0[p]_{\cO[\zeta_N]},\mu_{p,\cO[\zeta_N]}).
\end{split}
\end{equation}
When taking $\Qbar_p$-points, these give isomorphisms of
$\TT_{\FF_p,\mbar}$-modules, i.e.\ in particular of $\FF_p$-vector
spaces. We will from now on identify $\mu_{p,\cO[\zeta_N]}(\Qbar_p)$ with
$\FF_p$ and the group homomorphisms on~$\Qbar_p$-points above with
$\FF_p$-linear ones.

The final ingredient in the proof is the fact
that $G^e(\Qbar_p)[\m] = G^e[p](\Qbar_p)[\mbar]$ is a $1$-dimensional
$L := \TT_{\FF_p}/\mbar$-vector space by \cite{gross}, 
Prop.~12.8~(5) and~12.9~(4).
We now quotient the first isomorphism of Equation~\ref{grossvier}
(on $\Qbar_p$-points) by~$\mbar$
and obtain
$$ G^0[p](\Qbar_p)/\mbar \cong \Hom_{\FF_p}(G^e[p](\Qbar_p)[\mbar],\FF_p) 
\cong \Hom_{\FF_p}(L,\FF_p),$$
which is a $1$-dimensional $L$-vector space. Consequently, Nakayama's Lemma
applied to the finitely generated $\TT_{\FF_p,\mbar}$-module $G^0[p](\Qbar_p)$ yields
a surjection $\TT_{\FF_p,\mbar} \twoheadrightarrow G^0[p](\Qbar_p)$.
From \cite{kilford_wiese}, Prop~4.7, it follows that
$2 \dim_{\FF_p} \TT_{\FF_p,\mbar} = \dim_{\FF_p} H^1_{\text{par}} (\Gamma_1(Np),\FF_p)_\mbar$.
As we also have 
$$H^1_{\text{par}} (\Gamma_1(Np),\FF_p)_\mbar \cong J_1(Np)(\CC)[p]_\mbar \cong
G[p](\Qbar_p),$$
we obtain $\dim_{\FF_p} \TT_{\FF_p,\mbar} = \dim_{\FF_p} G^0[p](\Qbar_p)$ and,
thus, $\TT_{\FF_p,\mbar} \cong G^0[p](\Qbar_p)$, as desired.
\qed
\medskip

The following result together with very helpful hints on its proof
(i.e.\ the preceding proposition) was suggested by Kevin Buzzard.
See also the discussion before \cite{emerton}, Proposition~6.3,
and~\cite{mazur-appendix}.

\begin{cor}\label{gormult}
  Assume any of the two cases of Situation~\ref{sit} and let~$G$ be the
  $p$-divisible group attached to~$\m$. 
  Then there is the exact sequence
  $$ 0 \to \TT_{\FF_p,\mbar} 
       \to G[p](\Qbar) 
       \to \TT_{\FF_p,\mbar}^\vee 
       \to 0$$ 
  of $\TT_{\FF_p,\mbar}$-modules, where the dual is the $\FF_p$-linear dual.
\end{cor}

\pf
Substituting the isomorphism of Prop.~\ref{propq} into the second isomorphism 
of Equation~\ref{grossvier} (on $\Qbar_p$-points) gives
$$ G^e[p](\Qbar_p) \cong \Hom(\TT_{\FF_p,\mbar},\FF_p)$$
as $\TT_{\FF_p,\mbar}$-modules, 
whence the corollary follows from Equation~\ref{grossdrei}.
\qed
\medskip

The following proposition is similar in spirit to Proposition~\ref{propq}.
It will not be needed in the sequel.

\begin{prop}\label{tangent}
  Assume any of the two cases of Situation~\ref{sit} and let~$G$ be
  the $p$-divisible group attached to~$\m$. Then
  $G^0[p](\FF_p(\epsilon))$ and $\TT_{\FF_p, \mbar}$ are isomorphic as
  $\TT_{\FF_p,\mbar}$-modules.
\end{prop}

\pf 
We only give a sketch. Since $G^0[p](\FF_p)$ consists of the origin as
unique point, $G^0[p](\FF_p(\epsilon))$ coincides with the tangent
space at~$0$ of $G^0_{\FF_p}[p]$. The latter, however, is equal to the
tangent space at~$0$ of~$G_{\FF_p}[p]$. On the other hand, its dual,
the cotangent space at~$0$ of~$G_{\FF_p}[p]$, is isomorphic to
$S_p(\Gamma_1(N),\FF_p)_\mbar$. For Situation~(II) this is well-known.
In Situation~(I) we quote \cite{edixhoven-serre}, Eq.~6.7.1 and~6.7.2,
as well as \cite{gross}, Prop.~8.13 (note that the ordinarity
assumption kills the second summand in that proposition).
Consequently, $G^0[p](\FF_p(\epsilon))$ is isomorphic to the Hecke
algebra on $S_p(\Gamma_1(N),\FF_p)_\mbar$ as a Hecke module. In
\cite{kilford_wiese}, Prop.~2.3, it is shown that this algebra
is~$\TT_{\FF_p,\mbar}$.  
\qed 
\medskip

From Prop.~\ref{propq} and the reduction of points used in the direct
proof of Theorem~\ref{thmtp} we can also conclude an isomorphism
$G^0[p](\Fbar_p(\epsilon)) \cong \TT_{\FF_p,\mbar}$.

\section{Comparing Frobenius and the Hecke operator $T_p$}\label{frob}

The aim of this section is to prove that the Hecke operator~$T_p$ and
the Frobenius at~$p$ coincide on the unramified $\Qbar_p$-points
of~$G^0[p]$. 

\begin{thm}\label{thmtp}
  Assume any of the two cases of Situation~$\ref{sit}$ and let
  $G^0_\cO$ be the $p$-divisible group of Thm.~\ref{propgross}. The
  action of the geometric Frobenius on the points
  $G^0_\cO[p](\QQ_p^\nr(\zeta_p))$ is the same as the action of the
  Hecke operator~$T_p$.
\end{thm}

Using the Eichler-Shimura congruence relation in Situation~(II) and
the reduction of a well-known semi-stable model of the modular curve
in Situation~(I), the proof is quickly reduced to comparing the
geometric Frobenius and Verschiebung on the special fibre of~$G^0[p]$.
This comparison has been worked out conceptually by Niko Naumann
in Appendix~\ref{naumann} using Fontaine's theory of Honda systems in a
general setting. We also give a direct elementary proof. 
The idea of that proof is to work with the
tangent space at~$0$ over~$\Fbar_p$, in order to have an injective
reduction map from characteristic zero to the finite field. On the
special fibre elementary computations then suffice.
\medskip

{\bf Direct proof.}
We know that ${}^t G^0[p] = \Spec(A)$ is a finite \'etale group scheme
over~$\cO$ such that 
${}^t G^0[p] \times_\cO \ZZ_p^\nr[\zeta_p] \cong (\ZZ/p\ZZ)^h_{\ZZ^\nr_p[\zeta_p]}$,
i.e.\ $A \otimes_\cO \ZZ_p^\nr[\zeta_p] \overset{\alpha}{\cong} \prod \ZZ_p^\nr[\zeta_p]$.
If $p=2$, we put $\zeta_2 = -1$.
We obtain a reduction map
\begin{equation}\label{eqred}
  \Hom_{\grs/ \ZZ_p^\nr[\zeta_p]}
      ({}^t G^0[p] \times_\cO \ZZ_p^\nr[\zeta_p], \mu_{p,\ZZ_p^\nr[\zeta_p]}) \to 
  \Hom_{\grs/ \Fbar_p(\epsilon)} 
      ({}^t G^0[p] \times_{\cO} \Fbar_p(\epsilon), \mu_{p,\Fbar_p(\epsilon)})
\end{equation}
from the commutative diagram
$$ \xymatrix@=1.3cm{
\ZZ^\nr_p[\zeta_p][X]/(X^p-1)  \ar@{<-}[r]^(.45){\zeta_p \mapsfrom Y} \ar@{->}[d] &
\ZZ^\nr_p[X,Y]/(X^p-1,Y^p-1)   \ar@{->}[r]^(.55){Y \mapsto 1+\epsilon} \ar@{->}[d] &
\Fbar_p(\epsilon)[X]/(X^p-1) \ar@{->}[d] \\
\prod \ZZ^\nr_p[\zeta_p]   \ar@{<-}[r]^(.45){\zeta_p \mapsfrom Y} &
\prod \ZZ^\nr_p[Y]/(Y^p-1) \ar@{->}[r]^(.55){Y \mapsto 1+\epsilon} &
\prod \Fbar_p(\epsilon).}$$
Any morphism of group schemes 
${}^t G^0[p] \times_\cO \ZZ_p^\nr[\zeta_p] \to \mu_{p,\ZZ_p^\nr[\zeta_p]}$
corresponds to a Hopf algebra homomorphism as in the left column.
It is easy to see that it has a unique lifting to a homomorphism as
in the central column,
so that it gives a homomorphism in the right column.
Explicitly, a map in the left column is uniquely determined by
the image of~$X$, which is of the form
$(\zeta_p^{i_1},\dots,\zeta_p^{i_{hp}})$ for
some~$i_j \in \{0,\dots,p-1\}$. The corresponding map in the right column
sends~$X$ to $(1+i_1\epsilon,\dots,1+i_{hp}\epsilon)$.
Hence, the reduction map~\ref{eqred} is injective.
It is also compatible for the action induced by the Hecke correspondences.
In fact, for $p > 2$, one can pass directly from the left hand side column
to the right hand side via the map $\ZZ_p^\nr[\zeta_p] \to \Fbar_p(\epsilon)$,
sending $\zeta_p$ to~$1+\epsilon$.

Next, we describe the geometric Frobenius on the points ${}^t
G[p](\QQ_p^\nr(\zeta_p))$ and ${}^t G[p](\Fbar_p(\epsilon))$.
We consider the commutative diagram
$$ \xymatrix@=0.62cm{
\Hom_{\grs/\ZZ_p^\nr[\zeta_p]}({}^t G^0[p] \times \ZZ_p^\nr[\zeta_p],\mu_{p, \ZZ_p\nr[\zeta_p]})
 \ar@{->}[r]^(.7)\sim \ar@{->}[d]^\sim&
(A \otimes \ZZ_p^\nr[\zeta_p])^\gl \ar@{->}[r]^(.5)\sim \ar@{)->}[d] &
{}^tG[p](\ZZ_p^\nr[\zeta_p]) \ar@{)->}[d] \\
\Hom_{\ZZ_p^\nr[\zeta_p]-\Halg}(\ZZ_p^\nr[\zeta_p][X]/(X^p-1),A \otimes \ZZ_p^\nr[\zeta_p])
 \ar@{)->}[r] &
A \otimes \ZZ_p^\nr[\zeta_p] \ar@{->}[r]^(.4){\sim \; \ev}  &
\Hom_\cO({}^t A, \ZZ_p^\nr[\zeta_p]).}$$
It is well-known that a Hopf algebra homomorphism 
$\psi: \ZZ_p^\nr[\zeta_p][X]/(X^p-1) \to A \otimes_{\cO} \ZZ_p^\nr[\zeta_p]$
is uniquely given by the ``group-like element'' 
$\psi(X)= \sum a_i \otimes s_i$, giving the upper left bijection.
On the bottom right, we have the evaluation isomorphism
$A \otimes_\cO \ZZ_p^\nr[\zeta_p] \to \Hom_\cO(\Hom_\cO(A,\cO),\ZZ_p^\nr[\zeta_p])$
which is defined by $\ev(a \otimes s) (\varphi) = \varphi(a)s$.
We use that as $\cO$-modules ${}^t A = \Hom_\cO(A,\cO)$ 
with $G^0[p] = \Spec({}^t A)$, as well as the freeness of~$A$.
It is also well-known that the evaluation map gives rise to the upper right bijection.

Let now $\phi$ be the geometric Frobenius in $\Gal(\QQ_p^\nr(\zeta_p)/\QQ_p(\zeta_p))$.
Its action on $\Hom_\cO({}^t A, \ZZ_p^\nr[\zeta_p])$ is by composition.
Via the evaluation map it is clear that $\phi$ acts on an element
$a \otimes s \in A \otimes_\cO \ZZ_p^\nr[\zeta_p]$ by sending it
to $a \otimes \phi(s)$. Consequently, the morphism~$\psi^\phi$ on the left
which is obtained by applying~$\phi$ to~$\psi$ is uniquely determined by 
$\psi^\phi(X) = \sum a_i \otimes \phi(s_i)$.
A similar statement holds for the reduction. 
We note that this implies the compatibility of the
reduction map with the $\phi$-action.

Next we show that the action of geometric Frobenius on the tangent
space~$G^0[p](\Fbar_p(\epsilon))$ coincides with the action induced by
Verschiebung on~$G^0_{\FF_p}[p]$. The \'etale algebra~$A \otimes_\cO
\FF_p$ can be written as a product of algebras of the form $\FF_p[X]/(f)$ 
with $f$ an irreducible polynomial. 
An elementary calculation on the underlying rings gives
the commutativity of the diagram
\begin{equation}\label{diag}
 \xymatrix@=1.0cm{
\FF_p[X]/(f) \otimes_{\FF_p} \Fbar_p(\epsilon) \ar@{->}[r]^{F \otimes 1} \ar@{->}[d]^\alpha &
\FF_p[X]/(f) \otimes_{\FF_p} \Fbar_p(\epsilon) \ar@{->}[r]^{1 \otimes \phi^{-1}} &
\FF_p[X]/(f) \otimes_{\FF_p} \Fbar_p(\epsilon) \ar@{->}[d]^\alpha \\
\prod_{i=1}^d \Fbar_p(\epsilon) \ar@{->}[rr]^{\prod \phi^{-1}} &&
\prod_{i=1}^d \Fbar_p(\epsilon),}
\end{equation}
where $F$~denotes the absolute Frobenius on~${}^t G^0_{\FF_p}[p]$
(defined by $X \mapsto X^p$), which by duality gives the Verschiebung
on~$G^0_{\FF_p}[p]$. 
We point out that $\phi$ leaves~$\epsilon$ invariant.
Any $\Fbar_p(\epsilon)$-Hopf algebra homomorphism 
$\psi: \Fbar_p(\epsilon)[X]/(X^p-1) \to A \otimes_{\cO} \Fbar_p(\epsilon)$ 
is uniquely given by $\psi(X) = \sum_i a_i \otimes s_i$, and under
the identification 
$A \otimes_{\cO} \Fbar_p(\epsilon) \cong \prod_{j=1}^{hp} \Fbar_p(\epsilon)$ 
we get $\psi(X) = (1+i_1\epsilon,\dots,1+i_{hp}\epsilon)$, which is invariant under the
arithmetic Frobenius of the bottom row of~\ref{diag}. Hence,
$\phi^{-1}(F(\sum_i a_i \otimes s_i)) = \sum_i a_i \otimes s_i$, so
that $F(\sum_i a_i \otimes s_i) = \sum_i a_i \otimes \phi(s_i)$. This
proves that the geometric Frobenius and Verschiebung coincide.

We now finish the proof.
In Situation~(II) for $p=2$, the Eichler-Shimura relation $T_p =
\langle p \rangle F + V$ holds on the special fibre of $G[p]$
(see~\cite{gross}, proof of Prop.~12.8~(2)). Since~$F$ is zero on
$G^0_{\FF_p}[p]$, we get $T_p = V$ on it. As we have seen right above
that $V$ coincides with~$\phi$ on $G_{\FF_p}^0(\Fbar_p(\epsilon))$,
we obtain the theorem for~$p=2$.

In Situation~(I) we know that $G^0_{\FF_p}[p]$ is naturally part of
the $p$-torsion of the Jacobian of the Igusa curve $I_1(N)_{\FF_p}$;
but on the Igusa curve Verschiebung acts as~$T_p$ (see the proof of \cite{gross},
12.9~(2), for both these facts). Hence, we can argue as above and
get the theorem also for $p>2$.
\qed 
\medskip

{\bf More conceptual proof.}
In both situations, Theorem~\ref{thm} of Naumann gives an isomorphism
between $G^0[p](\QQ_p^\nr(\zeta_p))$ and the Dieudonn\'e module~$M$
attached to the special fibre~$G^0_{\FF_p}[p]$. Under this isomorphism
the geometric Frobenius $\phi \in
\Gal(\QQ_p^\nr(\zeta_p)/\QQ_p(\zeta_p))$ on
$G^0[p](\QQ_p^\nr(\zeta_p))$ is identified with Verschiebung on the
Dieudonn\'e module. The isomorphism is compatible with the Hecke
action.
Using the same citations as at the end of the direct proof
one immediately concludes that the equality $T_p = V$ holds 
on the Dieudonn\'e module~$M$, finishing
the proof.
\qed
\medskip

\begin{rem}
\begin{enumerate}[(a)]
\item Conceptually, taking $\ZZ_p^\nr[\zeta_p]$-points is the same as
  taking $\ZZ_p^\nr$-points of the Weil restriction from~$\cO$
  to~$\ZZ_p$.

\item For a representation~$\rho_\mbar$ which is unramified at~$p$ one
  knows that the arithmetic Frobenius~$\Frob_p$ satisfies $X^2 - T_p X
  + \langle p \rangle = 0$. This is in accordance with
  Theorem~\ref{thmtp}. For, it gives that $\Frob_p$ acts on
  $G^0[p](\Qbar_p)$ as~$a_p^{-1}$. Due to his conventions, Gross must
  still twist his representation by the determinant
  character~$\epsilon$, so that $\Frob_p$ acts as $\epsilon(p)/a_p$.
  This coincides with Deligne's description of the restriction
  of~$\rho_\mbar$ to a decomposition group at~$p$ (see, for instance,
  \cite{edixhoven-serre}, Thm.~2.5, or \cite{gross}, Prop.~12.1).
 
\end{enumerate}
\end{rem}

\section{Application to multiplicities}\label{conclude}

We first state a slight strengthening of a well-known theorem by
Boston, Lenstra and Ribet.

\begin{prop}[Boston, Lenstra, Ribet]\label{blr}
  Assume any of the two cases of Situation~\ref{sit}.
  Let~$m$ be an integer and $\FF = \TT_{\FF_p,\mbar}/\mbar$.
  Then the $\FF[\Gal(\Qbar/\QQ)]$-module
  $J(\Qbar)[\m^{(m)}]$ is the direct sum of $r$ copies of~$\rho_\mbar
  \otimes \epsilon^{-1}$ for some~$r \ge 1$ and Dirichlet character
  $\epsilon = \det(\rho_\mbar)$.
  
  The integer~$r$ is called the {\em multiplicity of $\rho_\mbar$
  on~$J(\Qbar)[\m^{(m)}]$}. If $m=1$, it is just called the 
  {\em multiplicity of $\rho_\mbar$}.
\end{prop}

\pf 
The same proof as in the original proposition works. More precisely,
one considers the two representations $\rho_\mbar: \Gal(\Qbar/\QQ) \to
\GL_2(\FF)$ and $\sigma: \Gal(\Qbar/\QQ) \to
\Aut(J(\Qbar)[\m^{(m)}])$. By Chebotarev's density theorem we know
that every conjugacy class of $\Gal(\Qbar/\QQ)/\ker(\sigma \otimes
\epsilon)$ is hit by a Frobenius element $\Frob_l$ for some $l \nmid Npm$.

The Eichler-Shimura congruence relation $T_l = \langle l \rangle F +
V$ holds on~$J_{\FF_l}$ (taking $J$ here over $\ZZ[\frac{1}{Np}]$) for
all primes $l \nmid Npm$.  Hence, the minimal polynomial of~$\Frob_l$
on the Jacobian divides $X^2 - T_l/\langle l \rangle \cdot X + l /
\langle l \rangle$.  But $T_l$ acts as~$a_l$ on $J(\Qbar)[\m^{(m)}]$
and $X^2 - a_l X + \epsilon(l) l$ (with $T_l \equiv a_l \mod \mbar$)
is the characteristic polynomial of
$\rho_\mbar(\Frob_l)$.  Consequently, $(\sigma \otimes \epsilon)(g) $ is
annihilated by the characteristic polynomial of~$\rho_\mbar(g)$ for
all $g\in \Gal(\Qbar/\QQ)$. Hence, Theorem~1 of~\cite{blr} gives the
result.
\qed
\medskip

The next corollary says that one can read off multiplicities from
properties of Hecke algebras.

\begin{cor}\label{multgor}
  Assume any of the two cases of Situation~\ref{sit}. Let $r$
  be the multiplicity of~$\rho_\mbar$.
  Then the relation
  $$ r = \frac{1}{2} 
    (\dim_\FF \TT_{\FF_p,\mbar}[\mbar] + 1)$$
  holds, where $\FF = \TT_{\FF_p,\mbar}/\mbar$.
\end{cor}

\pf
We note that in~\cite{buzzard_appendix} Buzzard explains
the exactness of the sequence
$$ 0 \to G^0(\Qbar_p)[\m] \to G(\Qbar_p)[\m] \to G^e(\Qbar_p)[\m] \to 0.$$
Via Corollary~\ref{gormult} we obtain the exact sequence
$$ 0 \to \TT_{\FF_p,\mbar} [\mbar] 
\to J_1(Np)(\Qbar_p)[\m] \to 
\big(\TT_{\FF_p,\mbar}/\mbar\big)^\vee \to 0,$$
from which one reads off the claim by counting dimensions.
\qed

\begin{thm}\label{mainthm}
  Assume any of the two cases of Situation~\ref{sit} and that
  $\rho_\mbar$ is of weight one.
  Then the following statements are equivalent.

\begin{enumerate}[(a)]
\item The representation $\rho_\mbar$ comes from a Katz cusp form of
  weight~$1$ on $\Gamma_1(N)$ over~$\FF_p$ and the multiplicity
  of~$\rho_\mbar$ is~$1$.

\item $\TT_{\FF_p,\mbar}[\mbar] \subsetneqq
  \TT_{\FF_p,\mbar}[\mbar^{(p)}]$

\item $T_p$ does not act as scalars on
  $\TT_{\FF_p,\mbar}[\mbar^{(p)}]$.

\item The multiplicity of~$\rho_\mbar$ is~$1$, its multiplicity on
  $J(\Qbar)[\mbar^{(p)}]$ is~$2$, and $\rho_\mbar(\Frob_p)$ is
  non-scalar.

\end{enumerate}

\end{thm}

\pf
$(a) \Rightarrow (b):$ By Cor.~\ref{multgor} and Nakayama's Lemma
$\TT_{\FF_p,\mbar}$ is Gorenstein, i.e.\ it is isomorphic to its
dual as a module over itself. By the $q$-expansion principle,
the dual is $S_p(\Gamma_1(N),\FF_p)_\mbar$.
By \cite{edixhoven}, Prop.~6.2, the existence of a corresponding
weight~$1$ form is equivalent to $S_p(\Gamma_1(N),\FF_p)_\mbar[\mbar^{(p)}]$
being $2$-dimensional. This establishes~$(b)$, since by the $q$-expansion
principle $S_p(\Gamma_1(N),\FF_p)_\mbar[\mbar]$ is $1$-dimensional.

$(b) \Rightarrow (c):$ This is evident.

$(c) \Rightarrow (d):$ First of all, $\TT_{\FF_p,\mbar}[\mbar^{(p)}]$ is
at least~$2$-dimensional (as $\TT_{\FF_p,\mbar}/\mbar$-vector space).
From Theorem~\ref{thmtp} we know that $T_p$ acts as the inverse of~$\Frob_p$
on $G^0[p](\Qbar)$.
We conclude that $\rho_\mbar(\Frob_p)$ cannot be scalar.
On $\TT_{\FF_p,\mbar}[\mbar]$ the action of $T_p$ is by scalars.
If the multiplicity~$r$ of~$\rho_\mbar$ were not~$1$,
then $\TT[\mbar] = G^0[p](\Qbar)[\mbar]$ would have dimension 
$2r-1>1$ (by the proof of Cor.~\ref{multgor}). 
From Theorem~\ref{thmtp} we obtain a contradiction.
We note that this argument, showing that $\rho_\mbar(\Frob_p)$ being
non-scalar implies that the multiplicity of~$\rho_\mbar$ is~$1$,
did not use the statement of~(c).
If the multiplicity~$s$ of~$\rho_\mbar$ on $J(\Qbar)[\mbar^{(p)}]$ were
bigger than~$2$, then $\TT_{\FF_p,\mbar}[\mbar^{(p)}]$ would be
at least $4$-dimensional. Then it follows that is must contain
at least two linearly independent eigenvectors for~$T_p$, contradicting
the fact that $\TT_{\FF_p,\mbar}[\mbar]$ is $1$-dimensional.

$(d) \Rightarrow (a):$ Clearly, $\mbar \neq \mbar^{(p)}$.
Hence, $\TT_{\FF_p,\mbar} / \mbar \neq \TT_{\FF_p,\mbar} / \mbar^{(p)}$
and, dually, 
$$S_p(\Gamma_1(N),\FF_p)_\mbar [\mbar] \subsetneqq 
  S_p(\Gamma_1(N),\FF_p)_\mbar [\mbar^{(p)}],$$
which implies the existence of a corresponding weight~$1$ form,
again by \cite{edixhoven}, Prop.~6.2.
\qed
\medskip

In~\cite{buzzard_appendix} Buzzard proved that the multiplicity of~$\rho_\mbar$
is~$1$ if $\rho_\mbar(\Frob_p)$ is non-scalar. We obtain that this
is in fact an equivalence (under a standard assumption in the case $p=2$).

\begin{cor}\label{coreins}
  Assume any of the two cases of Situation~\ref{sit} and that
  $\rho_\mbar$ is of weight one. If $p=2$, also assume that a weight~$1$
  Katz form of level~$N$ exists which gives rise to~$\rho_\mbar$.

  Then the multiplicity of~$\rho_\mbar$ is~$1$ if and
  only if $\rho_\mbar(\Frob_p)$ is non-scalar.
\end{cor}

\pf
By \cite{edixhoven-serre}, Theorem~4.5, together with the remark
at the end of the introduction to that article, the existence
of the corresponding weight~$1$ form is also guaranteed for~$p>2$.
If the multiplicity is~$1$, Theorem~\ref{mainthm} gives that
$\rho_\mbar(\Frob_p)$ is non-scalar. On the other hand, if
$\rho_\mbar(\Frob_p)$ is non-scalar, the argument 
used in the implication $(c) \Rightarrow (d)$ of 
Theorem~\ref{mainthm} shows that the multiplicity is~$1$.
\qed

\begin{cor}\label{corzwei}
  Assume any of the two cases of Situation~\ref{sit}. If $p=2$, also
  assume that if $\rho_\mbar$ is of weight one, then a weight~$1$ Katz
  form of level~$N$ exists which gives rise to~$\rho_\mbar$.

  Then the multiplicity of~$\rho_\mbar$ on $J[\m^{(p)}]$ is~$1$ if and
  only if $\rho_\mbar$ is ramified at~$p$.
\end{cor}

\pf
If $\rho_\mbar$ is ramified at~$p$, the result is Theorem~6.1
of~\cite{buzzard_appendix}.
Suppose now that $\rho_\mbar$ is unramified at~$p$. If
$\rho_\mbar(\Frob_p)$ is scalar, the corollary follows from
Corollary~\ref{coreins}. If $\rho_\mbar(\Frob_p)$ is non-scalar, 
then the result follows from Corollary~\ref{coreins} together with
the implication $(a) \Rightarrow (d)$ of Theorem~\ref{mainthm}.
\qed 

\begin{cor}\label{cordrei}
  Assume any of the two cases of Situation~\ref{sit} and that $\rho_\mbar$
  is of weight one.
  Assume also that the multiplicity of~$\rho_\mbar$ on
  $J(\Qbar)[\mbar^{(p)}]$ is~$2$. Then the following statements
  are equivalent.

\begin{enumerate}[(a)]
\item The multiplicity of~$\rho_\mbar$ is~$1$ and a weight~$1$ Katz
  form of level~$N$ exists which gives rise to~$\rho_\mbar$.
\item $\rho_\mbar(\Frob_p)$ is non-scalar.
\end{enumerate}
\end{cor}

\pf
We have seen the implication $(a) \Rightarrow (b)$ above. As in
the proof of Thm.~\ref{mainthm}, we obtain from $\rho_\mbar(\Frob_p)$
being non-scalar that the multiplicity of~$\rho_\mbar$ is~$1$. 
From the assumption the inequality
$\mbar \neq \mbar^{(p)}$ follows, implying the existence
of the weight~$1$ form as above by~\cite{edixhoven}, Prop.~6.2.
\qed
\medskip

If one could prove that the multiplicity of~$\rho_\mbar$ on $J(\Qbar)[\m^{(p)}]$ 
is always equal to~$2$ in the unramified situation, Corollary~\ref{cordrei}
would extend weight lowering for~$p=2$ to~$\rho_\mbar(\Frob_p)$
being non-scalar.

\newpage
\appendix

\newcommand{\Z}{\mathbb{Z}}
\newcommand{\Qp}{\mathbb{Q}_p}
\newcommand{\F}{{\mathcal F}}
\newcommand{\Fp}{\mathbb{F}_p}
\newcommand{\Fpb}{\overline{\mathbb{F}}_p}
\newcommand{\Zp}{\mathbb{Z}_p}
\newcommand{\Gl}{\mathrm{Gl}}
\newcommand{\mm}{\mathrm{m}}
\newcommand{\M}{{\mathcal M}}
\newcommand{\V}{{\mathcal V}}
\newcommand{\N}{\mathrm{N}}
\renewcommand{\O}{{\mathcal O}}
\newcommand{\D}{\mathrm{D}}
\newcommand{\CW}{\mathrm{CW}}
\newcommand{\CWF}{CW_{\Fp}(\pi{\mathcal O}/\pi^2{\mathcal O})}

\section{Appendix}\label{naumann}
\begin{large}\begin{center}
By Niko Naumann\footnote{NWF 1-Mathematik,
Universität Regensburg,
D-93040 Regensburg,
Germany,\\ \hspace*{.5cm} niko.naumann@mathematik.uni-regensburg.de}
\end{center}\end{large}

Let $p$ be a prime, $A:=\Zp$, $A':=\Zp[\zeta_p]$, $K:=\Qp$,
$K':=\Qp(\zeta_p)$ and $K'\subseteq\overline{K}$ an algebraic closure.
We have the inertia sub-group $I\subseteq
G_{K'}:=\Gal(\overline{K}/K')$ and for a $G_{K'}$-module $V$ we denote
by $\tau$ the geometric Frobenius acting on the inertia invariants
$V^I$. If $G/A'$ is a finite flat group-scheme, always assumed to be
commutative, we denote by $M$ the Dieudonn\'e-module of its special
fiber and by $V:M\to M$ the Verschiebung.

\begin{thm}\label{thm} Let $G/A'$ be a finite flat group-scheme which is connected 
  with \'etale Cartier-dual and annihilated by multiplication with
  $p$. Then $G(\overline{K})^I=G(\overline{K})$ and there is an
  isomorphism $\phi:G(\overline{K})^I\to M$ of $\Fp$-vector spaces
  such that $\phi\circ\tau=V\circ\phi$.
\end{thm}

The assumption that $pG=0$ cannot be dropped in Theorem \ref{thm}:

\begin{prop}\label{prop} For every $n\ge 2$ there is a finite flat group-
  scheme $G/A'$ of order $p^n$ which is connected with an \'etale dual
  and such that $G(\overline{K})^I\simeq\Z/p\Z$ with $\tau$ acting
  trivially and $V\neq 1$ on the Dieudonn\'e-module of the special
  fiber of $G$.
\end{prop}

{\bf Proof of Theorem \ref{thm}.} Denoting by $G'$ the Cartier-dual of
$G/A'$ we have an isomorphism of $G_{K'}$-modules
\[
G(\overline{K})\simeq\Hom(G'(\overline{K}),\mu_{p^{\infty}}(\overline{K}))\stackrel{(pG'=0)}{=}\Hom(G'(\overline{K}),\mu_{p}(\overline{K})).\]

Since $G'(\overline{K})$ is unramified because $G/A'$ is \'etale and
$\mu_{p}(\overline{K})$ is unramified because $\zeta_p\in K'$ we see
that $G(\overline{K})^I=G(\overline{K})$. Letting $p^n$ denote the
order of $G$ we have
\[
\dim_{\Fp}(G(\overline{K})^I)=\dim_{\Fp}(G(\overline{K}))=n=\dim_{\Fp}(M).\]

In the rest of the proof we use the explicit quasi-inverse to J.-M.
Fontaine's functor associating with $G$ a finite Honda system in order
to determine
the action of $\tau$ on $G(\overline{K})^I$ \cite{Fo},\cite{Co}.

Let $(M,L)$ be the finite Honda system over $A'$ associated with
$G/A'$.  Recall that $M$ is the Dieudonn\'e-module of the special
fiber of $G$ and $L\subseteq M_{A'}$ is an $A'$-sub-module where
$M_{A'}$ is an $A'$-module
functorially associated with $M$ \cite[Ch. IV,\S 2]{Fo}.

We claim that $L=M_{A'}$: Let $\mm\subseteq A'$ denote the maximal
ideal. Using the notation of \cite[Section 2]{Co}, the defining
epimorphism of $A'$-modules $M_{A'}\to\coker(\F_M)$ factors through an
epimorphism $M_{A'}/\mm M_{A'}\to\coker(\F_M)$ because $\mm\cdot
\coker(\F_M)=0$ \cite[Lemma 2.4]{Co}. Denoting by $l$ the length of a
module we have
\[ l_{A'}(\coker(\F_M))\stackrel{\cite[2.4]{Co}}{=}l_A(\ker\, F)=
l_A(\ker(p:M\to M))=l_A(M)=n\]
because $\ker\, F=\ker(p:M\to M)$ since $V$ is bijective, and $pM=0$.
On the other hand, the canonical morphism of $A'$-modules
$\iota_M:M\otimes_A A'\to M_{A'}$ is an isomorphism by \cite[Ch. IV,
Proposition 2.5]{Fo} using again that $V$ is bijective. Thus
\[ l_{A'}(M_{A'}/\mm M_{A'})=l_{A'}(M\otimes_A A'/
\mm)=l_{A'}(M/pM)=l_A(M)=n \] and $M_{A'}/\mm
M_{A'}\stackrel{\simeq}{\to}\coker(\F_M)$. Since $L/\mm
L\stackrel{\simeq}{\to}\coker(\F_M)$ holds for every finite Honda
system we see that the inclusion $L\subseteq M_{A'}$
induces an isomorphism $L/\mm L\stackrel{\simeq}{\to}M_{A'}/\mm M_{A'}$
and Nakayama's lemma implies that $L=M_{A'}$.

Fix $\pi\in\overline{K}$ with $\pi^{p-1}=-p$, then $K'=K(\pi)$: This
is obvious for $p=2$ and for $p\neq 2$ it follows from local class
field theory and the norm computations
$\N^{K'}_K(\zeta_p-1)=\N^{K(\pi)}_K(\pi)=p$. Note that $\pi\in A'$ is
a local uniformizer. Let $K'^{ur}$ denote the completion of the
maximal unramified extension of $K'$ inside $\overline{K}$ and
$\O\subseteq K'^{ur}$
its ring of integers.

By \cite[Remarque on p. 218]{Fo} and the fact that $L=M_{A'}$ we see that 
reduction induces an isomorphism
\begin{equation}\label{red1}
  G(\overline{K})^I=G(K'^{ur})=G(\O)\stackrel{\simeq}{\to}
  \left\{ \phi\in\Hom_{\D_{\Fp}}(M,\CWF)\, \, |\, w'^c\circ\phi_{A'}=0 \right\} 
\end{equation}
where $\D_{\Fp}=\Fp[F,V]$ is the Dieudonn\'e-ring, $\CW$ denote
Witt-covectors \cite[Ch. II,\S 1]{Fo},
\[ w'^c:\CWF_{A'}\to K'^{ur}/\pi^2\O\] is as in \cite[Ch. IV, \S
3]{Fo} and $\phi_{A'}:M_{A'}\to\CWF_{A'}$ is induced by $\phi$. By
construction of $w'^c$ we have, for every $\phi\in
\Hom_{\D_{\Fp}}(M,\CWF)$, a commutative diagram
\[\xymatrix{ M_{A'}\ar[r]^-{\phi_{A'}} & \CWF_{A'} \ar[r]^-{w'^c} & K'^{ur}/\pi^2\O \\
  & \CWF\otimes_A A'\ar@{->>}[u]^{\iota_{\CWF}} \ar[ur]^{\tilde{w}} & \CWF \ar[l] \ar[u]^{w^c} \\
  M\otimes_A A'\ar[uu]^{\iota_M}_{\simeq} \ar[ru]^-{\phi\otimes 1} & M
  \ar[ur]^-{\phi}\ar[l]} \]
in which
$w^c((x_{-n})_{n\ge0})=\sum\limits_{n=0}^{\infty}p^{-n}\hat{x}_{-n}^{p^n}$
with $\hat{x}_{-n}\in\pi\O$ lifting $x_{-n}$, $\tilde{w}=w^c\otimes 1$
is the $A'$-linear extension of $w^c$ and $\iota_{\CWF}$ is surjective
by \cite[Ch. IV, Proposition 2.5]{Fo} since $\CWF$ is $V$-divisible.
It is easy to see that we have
\begin{equation}\label{cl1} w'^c\circ \phi_{A'}=0\Leftrightarrow w^c\circ\phi=0.
\end{equation}

Combining (\ref{cl1}) and (\ref{red1}) we obtain an isomorphism 
\begin{equation}\label{red2}
G(\overline{K})^I\stackrel{\simeq}{\to}
\{ \phi\in\Hom_{\D_{\Fp}}(M,\CWF)\, | \, w^c\circ\phi=0 \}.
\end{equation}

Now we need to study $\ker(w^c)$. We will use the isomorphism of $\Fpb$- vector spaces
\begin{equation}\label{id}
\pi\O/\pi^2\O\stackrel{ : \pi}{\to}\O/\pi\O\simeq\Fpb
\end{equation}
to describe elements of $\CWF$ as covectors $(y_{-n})_{n\ge 0}$ with
$y_{-n}\in\Fpb$. Of course, since (\ref{id}) is not multiplicative,
some care has to be taken with this. We denote by
$\sigma:\Fpb\to\Fpb\, , \, \sigma(x)=x^p$ the absolute Frobenius and
claim that
\begin{equation}\label{ker}
\ker(w^c)=\{ (y_{-n})_n\, | \, y_{-n}\in\Fpb\, , \, y_{-1}=y_0^{\sigma^{-1}} \}.
\end{equation}

To see this, let $(x_{-n})_n\in\CWF$ be given, choose $\hat{x}_{-n}\in\pi\O$
lifting $x_{-n}$ and write $\hat{x}_{-n}=\pi \hat{y}_{-n}$ with $\hat{y}_{-n}\in\O$.
Then we compute in $K'^{ur}/\pi^2\O$:
\[ w^c((x_{-n}))=\sum\limits_{n=0}^{\infty}p^{-n}(\pi
\hat{y}_{-n})^{p^n}\stackrel{(\pi^{p-1}=-p)}{=}
\sum\limits_{n=0}^{\infty}(-1)^n\pi^{p^n-n(p-1)}\hat{y}_{-n}^{p^n}=\pi(\hat{y}_0-\hat{y}_{-1}^p),
\]
using that $p^n-n(p-1)\ge 2$ for all $n\ge 2$. Now (\ref{ker}) is obvious.

Next, we claim that the subset
\begin{equation}\label{subM}
\CWF\supseteq\M:=\{ (y_0^{\sigma^{-n}})_{n\ge 0}\, | \, y_0\in\Fpb \}
\end{equation}
is a $\D_{\Fp}$-sub-module. First note that $F=0$ on $\CWF$ so we will
consider it as a $\D_{\Fp}/F=\Fp[V]$-module in the following. Since
all products in $\pi\O/\pi^2\O$ are zero we have
$$(x_{-n})+(y_{-n})=(x_{-n}+y_{-n})$$ in $\CWF$ and
$\M$ is indeed a $\Fp$-sub-module, visibly stable under $V$.

We claim that the inclusion (\ref{subM}) induces an isomorphism
\begin{equation}\label{id2}
\Hom_{\Fp[V]}(M,\M)\stackrel{\simeq}{\to}
\{ \phi\in\Hom_{\D_{F_p}}(M,\CWF)\, | \, w^c\circ\phi=0 \}.
\end{equation}

Since $\M\subseteq\ker(w^c)$ by (\ref{ker}) we only need to see that a
$\Fp[V]$-linear morphism $$\phi:M\to\CWF$$ with
$\phi(M)\subseteq\ker(w^c)$ factors through $\M$: For every $m\in M$
and $n\ge 0$ we have, writing $\phi(m)=:(y_{-n})$ with
$y_{-n}\in\Fpb$,
\[ 0=w^c(\phi(V^nm))=w^c(V^n(\phi(m)))=w^c((\ldots,y_{-n-1},y_{-n})),
\]
thus $y_{-n-1}=y_{-n}^{\sigma^{-1}}$ by (\ref{ker}) and as this is true
for every $n\ge 0$ we get $\phi(m)\in\M$.

To proceed, note that
\begin{equation}\label{id3}
\M\to\Fpb\, , \, (y_0^{\sigma^{-n}})\mapsto y_0
\end{equation}
is an isomorphism of $\Fp[V]$-modules if one defines
$V(\alpha):=\alpha^{\sigma^{-1}}$ for $\alpha\in\Fpb$. Denoting by
$\Phi:G(\overline{K})^I\stackrel{\simeq}{\to}\Hom_{\Fp[V]}(M,\Fpb)$
the isomorphism obtained by combining (\ref{red2}), (\ref{id2}) and
(\ref{id3}), by construction we have a commutative diagram
\begin{equation}\label{final}
  \xymatrix{ G(\overline{K})^I \ar[r]^-{\Phi} \ar[d]^{\tau} & 
  \Hom_{\Fp[V]}(M,\Fpb) \ar[d]^{\Hom(V,\Fpb)} \\ 
    G(\overline{K})^I \ar[r]^-{\Phi} & \Hom_{\Fp[V]}(M,\Fpb).}
\end{equation}

Let $e_i$ (resp. $\phi_i$) ($1\leq i\leq n$) be an $\Fp$-basis of $M$
(resp. $\Hom_{\Fp[V]}(M,\Fpb)$) and define $Ve_i=:\sum_ja_{ij}e_j$,
hence $A:=(a_{ij})\in\Gl_n(\Fp)$,
$\psi_i:=\Hom(V,\Fpb)(\phi_i)=:\sum_j b_{ij}\phi_j$, hence
$B:=(b_{ij})\in\Gl_n(\Fp)$ and $C:=(\phi_i(e_j))\in\Gl_n(\Fpb)$. By
definition, $A$ is a representing matrix of $V:M\to M$ and by
(\ref{final}) $B$ is a representing matrix for $\tau$.
So we will be done if we can show that $A$ and $B$ are conjugate over $\Fp$.

From the computation
\[ \psi_i(e_j)=\phi_i(Ve_j)=\sum_k a_{jk}\phi_i(e_k)=\sum_k
b_{ik}\phi_k(e_j)\] we obtain $^tA=C^{-1}BC$. Now recall that over
every field $\kappa$ two square matrices with coefficients in $\kappa$
which are conjugate over an algebraic closure of $\kappa$ are
conjugate over $\kappa$ and, furthermore, that every square matrix
with coefficient in $\kappa$ is conjugate, over $\kappa$, to its
transpose. Hence $A$ is indeed conjugate to $B$ over $\Fp$.
\qed

\begin{rem}\label{rem}
  Inspecting the above proof we see that for $G/A'$ connected with
  \'etale dual (not necessarily annihilated by $p$) we have a
  commutative diagram
\[ \xymatrix{ G(\overline{K})^I\ar[r]^-{\Phi}_-{\simeq} \ar[d]^{\tau} & \Hom_{\Fp[V]}(M/FM,\Fpb)
\ar[d]^-{\Hom(V,\Fpb)} \\ 
G(\overline{K})^I \ar[r]^-{\Phi}_-{\simeq} & \Hom_{\Fp[V]}(M/FM,\Fpb).} \]
\end{rem}

{\bf Proof of Proposition \ref{prop}.} Define a finite Honda system over $A'$
by 
$$M:=\Z/p^n\Z, \; 1\neq V\in 1+p(\Z/p^n\Z)\subseteq(\Z/p^n\Z)^*=\Aut_{\Zp}(M),\; 
F:=pV^{-1},\; L:=M_{A'}.$$
It is easy to see that this is indeed a
finite Honda system. For the corresponding group $G/A'$ we have by
Remark \ref{rem}
\[ G(\overline{K})^I\simeq \Hom_{\Fp[V]}(M/FM,\Fpb)=\Fpb^{V=1}=\Fp\]
with trivial geometric Frobenius, note that $V$ is the identity on $M/FM$, but
$V\neq 1$.
\qed

\bibliographystyle{plain}
\bibliography{mult}

\end{document}